\documentclass{amsart}
\usepackage{amscd,amssymb,enumerate,amsmath,comment}
\usepackage[all]{xy}

\newcommand{\Hom}{\operatorname{Hom}\nolimits}

\newcommand{\MaxSpec}{\operatorname{MaxSpec}\nolimits}
\newcommand{\gr}{\operatorname{gr}\nolimits}
\newcommand{\HH}{\operatorname{HH}\nolimits}
\newcommand{\mo}{\mathfrak{o}}
\newcommand{\mt}{\mathfrak{t}}
\newcommand{\mm}{\mathfrak{m}}
\renewcommand{\L}{\Lambda}

\newcommand{\kar}{\operatorname{char}\nolimits}

\newcommand{\Ext}{\operatorname{Ext}\nolimits}
\newcommand{\N}{\mathcal{N}}
\newcommand{\R}{\mathcal{R}}

\newcommand{\rad}{\mathfrak{r}}

\newcommand{\Ann}{\operatorname{Ann}\nolimits}

\newtheorem{lem}{Lemma}[section]
\newtheorem{prop}[lem]{Proposition}

\newtheorem{thm}[lem]{Theorem}
\theoremstyle{definition}
\newtheorem{defin}[lem]{Definition}

\newtheorem{example}[lem]{Example}

\title[Support varieties for modules \dots]
{Support varieties for modules over stacked monomial algebras}
\thanks{The work in this paper was done while the first author was 
visiting the University of Leicester in 2008/2009. The first 
author would like to express his gratitude for the hospitality 
and support from the University of Leicester. The second author wishes to thank
the University of Leicester for granting study leave.}
\author[Furuya]{Takahiko Furuya}
\author[Snashall]{Nicole Snashall}
\address{Takahiko Furuya\\ Department of Mathematics\\
University of Leicester\\
University Road\\
Leicester, LE1 7RH\\
England}
\email{tf47@mcs.le.ac.uk}
\address{Nicole Snashall\\ Department of Mathematics\\
University of Leicester\\
University Road\\
Leicester, LE1 7RH\\
England}
\email{N.Snashall@mcs.le.ac.uk}



\begin{document}

\begin{abstract}
Let $\L$ be a finite-dimensional $(D,A)$-stacked monomial algebra.
In this paper, we give necessary and sufficient conditions for the
variety of a simple $\L$-module to be nontrivial. This is then used to give
structural information on the algebra $\L$, as it is shown
that if the variety of every simple module is nontrivial, then
$\L$ is a $D$-Koszul monomial algebra. We also provide examples of 
$(D,A)$-stacked monomial algebras which are not self-injective but
nevertheless satisfy the finite generation conditions {\bf (Fg1)} and 
{\bf (Fg2)} of \cite{EHSST}, from which we can characterize all modules with
trivial variety.
\end{abstract}

\maketitle

\section*{Introduction}\label{intro}

Let $K$ be an algebraically closed field, and let
$\L=K\mathcal{Q}/I$ be an indecomposable finite-dimensional algebra, where
${\mathcal Q}$ is a finite quiver and $I$ is an admissible ideal. Let $\L^e$ be 
the enveloping algebra $\L\otimes_K\L^{\rm op}$. The Hochschild cohomology
ring of $\L$ is defined to be $\HH^*(\L) = \Ext^*_{\L^e}(\L,\L) = \oplus_{n
\geq 0}\Ext^n_{\L^e}(\L, \L)$ with the
Yoneda product. We let $\N$ denote the ideal in $\HH^*(\Lambda)$ which is
generated by the homogeneous nilpotent elements.

In this paper we consider support varieties for modules over $(D,A)$-stacked
monomial algebras. This class was introduced by Green and
Snashall in \cite{GS}, and the $(D,A)$-stacked monomial algebras of infinite global
dimension were characterized in \cite{GS2} as precisely the monomial algebras 
for which every projective module in a minimal projective resolution of
$\L/\rad$ over $\L$ is generated in a single degree and for which the Ext
algebra is finitely generated as a $K$-algebra (but not
finite-dimensional). Note that we write $\rad$ for the 
Jacobson radical of $\L$. The class of $(D,A)$-stacked monomial algebras 
includes the Koszul monomial
algebras (equivalently, the quadratic monomial algebras) and the
$D$-Koszul algebras of Berger \cite{B}.

Support varieties for modules over any finite-dimensional algebra $\L$ were 
introduced by Snashall and Solberg in
\cite{SS}, where the support variety $V(M)$ of a finitely 
generated $\L$-module $M$ was defined by
\[V(M)=\{ \mm\in \MaxSpec\HH^*(\L)/\N\mid
\Ann_{\HH^*(\L)}\Ext^*_\L(M,M) \subseteq \mm'\}\] where $\mm'$
denotes the inverse image of $\mm$ in $\HH^*(\L)$. 
There is a unique maximal graded ideal $\mm_{\gr}$ in $\HH^*(\L)/\N$, 
and $\{\mm_{\gr}\} \subseteq V(M)$ for all non-zero finitely generated 
$\L$-modules $M$ (\cite[Proposition 3.4]{SS}). The variety of $M$ is 
then said to be trivial if $V(M) = \{\mm_{\gr}\}$.

It is clear that if $\L$ is of finite global dimension then every
homogeneous element of $\HH^*(\L)$ of strictly positive degree is
nilpotent, and so $\HH^*(\L)/\N \cong K$. It then follows that every module
has trivial variety. Since we are interested in the geometric and algebraic 
connections between varieties and homological properties of modules, we
restrict ourselves to algebras of infinite global dimension.

In this paper we study support varieties for modules over a $(D,A)$-stacked
monomial algebra $\L$. In Theorem \ref{simple_module_variety}, we give 
necessary and sufficient conditions for a simple $\L$-module to have trivial 
variety. This leads to Theorem \ref{all_simples_nontrivial}, where we provide new
structural information on the algebra $\L$, namely, we show that if every
simple $\L$-module has nontrivial variety then $A=1$ and so $\L$ is a $D$-Koszul
algebra. It was shown in \cite{EHSST}, that many of the properties of support
varieties for group algebras have analogues in the more general case under
the assumption that certain finiteness conditions {\bf (Fg1)} and {\bf (Fg2)} 
hold.
In particular, \cite[Theorem 2.5]{EHSST} shows that under these conditions, the algebra
is necessarily Gorenstein (that is, it has finite injective dimension as
both a left and a right module), and moreover, the variety of a module is 
trivial if and only if the module has finite projective dimension.  In Section
\ref{Section3}, we
give examples of $(D,A)$-stacked monomial algebras which are not 
self-injective algebras but nevertheless where conditions {\bf (Fg1)}
and {\bf (Fg2)} hold. This gives some new and interesting examples where the 
properties of support varieties for group algebras have analogues for algebras 
which are not self-injective.

\section{Background}

We fix the notation of this paper. 
A finite-dimensional algebra $\L = K\mathcal{Q}/I$ is a monomial algebra 
if $I$ is generated by a set of paths in $K{\mathcal Q}$ each of length at
least 2. We fix a minimal generating set $\rho$ for the ideal $I$, and refer to 
an element of $\rho$ as a {\it relation}. We assume further that $K$ is an
algebraically closed field with $\kar K\neq2$, and that $\L$ is
indecomposable of infinite 
global dimension.

An arrow $\alpha$ starts at the vertex $\mathfrak{o}(\alpha)$ and
ends at the vertex $\mathfrak{t}(\alpha)$; arrows in a path are read
from left to right. If $p = \alpha_1\alpha_2\cdots \alpha_n$ is a
path in $K{\mathcal Q}$ where $\alpha_1, \ldots , \alpha_n$ are arrows, 
then we set $\mo(p) = \mathfrak{o}(\alpha_1)$ and $\mt(p) =
\mathfrak{t}(\alpha_n)$. We denote the {\it length}
of a path $p$ by $\ell(p)$. An arrow $\alpha$ is the {\it first arrow} (respectively
{\it last arrow}) of a relation $r$ in $\rho$ if $r = \alpha p$
(respectively $r = p\alpha$) for some path $p$.

The definition of a $(D,A)$-stacked monomial algebra was given in \cite{GS} 
in terms of sets ${\mathcal R}^n$ which are constructed from  
overlaps of paths in ${\mathcal Q}$. These sets ${\mathcal R}^n$ are
also used to define a minimal projective bimodule resolution of $\L$. For
more details on overlaps see \cite{GS}.

\begin{defin}\cite[Definition 3.1]{GS}
Let $\L = K\mathcal{Q}/I$ be a finite-dimensional monomial algebra,
where $I$ has a minimal set of generators $\rho$. Then $\L$ is said to be 
a \emph{$(D,A)$-stacked monomial algebra} if there is some $D \geq 2$ and 
$A \geq 1$ such that, for all $n \geq 2$ and $R^n \in \R^n$,
$$\ell(R^n) = \left \{ \begin{array}{ll}
\frac{n}{2}D & \mbox{$n$ even}\\\\
\frac{(n-1)}{2}D + A & \mbox{$n$ odd.}
\end{array} \right.$$
In particular all relations in $\rho$ are of length $D$.
\end{defin}

The structure of the Hochschild cohomology ring modulo nilpotence of a
$(D,A)$-stacked monomial algebra was determined in \cite[Theorem 3.4]{GS}.
We use this result heavily in this paper, and therefore we now recall the 
notation and the result.

A \emph{closed path} $C$ in ${\mathcal Q}$ at the vertex $v$ is a 
non-trivial path $C$ in $K\mathcal{Q}$ with $C = vCv$ for some vertex $v$. 
If $C$ is a closed path at the vertex $v$ then we say that $v$ is
\emph{not internal} to $C$ if $C = v\sigma_1 v\sigma_2 v$ for paths
$\sigma_1$, $\sigma_2$ implies that $\sigma_1 = v$ or $\sigma_2 =
v$.

For $A \geq 1$, a {\it closed $A$-trail} $T$ in $\mathcal{Q}$ is a non-trivial
closed path $T = \alpha_0\alpha_1 \cdots \alpha_{m-1}$ in $K{\mathcal Q}$ 
such that $\alpha_0, \ldots , \alpha_{m-1}$ are all distinct paths of length $A$. 
By setting $T_0 = T$ and $$\begin{array}{rcl}
T_1 & = & \alpha_1 \cdots \alpha_{m-1}\alpha_0\\
T_2 & = & \alpha_2 \cdots \alpha_0\alpha_1\\
\vdots & & \\
T_{m-1} & = & \alpha_{m-1}\alpha_0 \cdots \alpha_{m-2},
\end{array}$$
then we say that $\{T_0, T_1, \ldots , T_{m-1}\}$
is a \emph{complete set of closed $A$-trails} on the $A$-trail $T = 
\alpha_0 \alpha_1 \cdots \alpha_{m-1}$.

Now, let $d \geq 2$ and write $d = Nm + l$ where $0 \leq l \leq m-1$ and
$N \geq 0$. For $t \in {\mathbb N}$, let $[t] \in \{0, 1, \ldots ,
m-1\}$ denote the residue of $t$ modulo $m$. Let $W =
T_0^N\alpha_0\alpha_1 \cdots \alpha_{l-1}$ with the conventions that
if $N = 0$ then $T^N_0 = \mo(\alpha_0)$ and if $l = 0$ then $W = T_0^N$. More
generally, for $k = 0, 1, \ldots , m-1$, define $\sigma^k(W) =
T_k^N\alpha_k\alpha_{k+1} \cdots \alpha_{k+l-1}$ with the
conventions that 
\begin{enumerate}
\item[(i)] if $t \geq m$ then $\alpha_t = \alpha_{[t]}$,
\item[(ii)] if $N = 0$ then $T^N_k = e_k$, and 
\item[(iii)] if $l = 0$ then $\sigma^k(W) = T_k^N$. 
\end{enumerate}
Define $\rho_T$ to be
the set
$$\rho_T = \{W, \sigma(W), \ldots , \sigma^{m-1}(W)\}.$$
We say that $\rho_T$ is the \emph{set of paths of length $dA$ that
are associated to the $A$-trail} $T$. Note that $\{W, \sigma(W),
\ldots , \sigma^{m-1}(W)\}$ is also the set of paths of length $dA$
that is associated to each $A$-trail $T_k$ for $k = 0, \ldots ,
m-1$.

\bigskip

Let $\L$ be a $(D,A)$-stacked monomial algebra of infinite global dimension.
Then \cite[Proposition 3.3]{GS} tells us that $D = dA$ for some $d \geq 2$.

Let $C_1, \ldots , C_u$ be all the closed paths in the quiver
$\mathcal{Q}$ at the vertices $v_1, \ldots , v_u$ respectively, such
that for each $C_i$ with $1\leq i \leq u$, we have $C_i\neq
p_i^{r_i}$ for any path $p_i$ with $r_i \geq 2$, $C_i^d \in \rho$, and 
there are no overlaps of $C_i^d$ with any relation in $\rho\setminus\{C_i^d\}$.
(Note that it follows that $\ell(C_i) = A$.)

Let $T_{u+1}, \ldots , T_r$ be all the distinct closed $A$-trails in
the quiver $\mathcal{Q}$ such that for each $T_i$ with $u+1 \leq i
\leq r$, the set $\rho_{T_i}$ of paths of length $D = dA$ which are
associated to the trail $T_i$ is contained in $\rho$ but, if $T_i =
\alpha_{i\, 0}\alpha_{i\, 1}\cdots\alpha_{i\, m_i-1}$, then each
path $\alpha_{i\, j}$ of length $A$ has no overlaps with any
relation in $\rho\setminus\rho_{T_i}$. (We assume that there is no 
repetition amongst these closed paths and
closed $A$-trails, that is, $\{C_1, \ldots , C_u\}\cap\{T_{u+1}, \ldots ,
T_r\} = \emptyset$.)

Then, \cite[Theorem 3.4]{GS} proves that 
$$\HH^*(\L)/\N \cong
K[x_1, \ldots , x_r]/\langle x_ax_b \mid a
\neq b\rangle$$ where
\begin{enumerate}
\item[(a)] for $j = 1, \ldots , u$, the vertices $v_1, \ldots , v_u$ are
distinct, and the element $x_j$ corresponding to the closed path
$C_j$ is in degree 2 and is represented by the map $P^2 \rightarrow
\L$ where, for $R^2\in\R^2$,
$$\mo(R^2)\otimes \mt(R^2) \mapsto \left \{
\begin{array}{ll} v_j & \mbox{if $R^2 = C_j^d$}\\
                   0 & \mbox{otherwise,}
\end{array} \right.$$
\item[(b)] and for $j = u+1, \ldots , r$, let $T_{j,0}, \ldots ,
T_{j,m_j-1}$ denote the complete set of $A$-trails on the closed
path $T_j$. Then the element $x_j$ corresponding to the closed
$A$-trail $T_j$ is, in the above notation, in degree $2\mu_j$ where
$\mu_j = m_j/\gcd(d, m_j)$ and is represented by the map $P^{2\mu_j}
\rightarrow \L$ where, for $R^{2\mu_j}\in\R^{2\mu_j}$,
$$\mo(R^{2\mu_j})\otimes \mt(R^{2\mu_j}) \mapsto \left \{
\begin{array}{ll} \mo(T_{j,k}) & \mbox{if $R^{2\mu_j} = T_{j,k}^{d/\gcd(d,
m_j)}$}\\
&\ \ \ \ \mbox{for $k = 0, \ldots , m_j-1$}\\
                   0 & \mbox{otherwise.}
\end{array} \right.$$
\end{enumerate}

Throughout this paper, whenever we refer to a closed path $C_i$ in the 
quiver $\mathcal{Q}$, for some $1 \leq i \leq u$, then we are assuming 
without further
comment that $C_i$ is a closed path at the vertex $v_i$, such that
$C_i\neq p_i^{r_i}$ for some path $p_i$ with $r_i \geq 2$, $C_i^d
\in \rho$ where $d = D/A$, and there are no overlaps of $C_i^d$ with
any relation in $\rho\setminus\{C_i^d\}$. Similarly, when we refer
to a closed $A$-trail $T_i$ in the quiver $\mathcal{Q}$, for some
$u+1 \leq i \leq r$, then we are assuming without further comment
that the set $\rho_{T_i}$ of paths of length $D$ which are
associated to the trail $T_i$ is contained in $\rho$ but, if $T_i =
\alpha_{i\, 0}\alpha_{i\, 1}\cdots\alpha_{i\, m_i-1}$, then each
path $\alpha_{i\, j}$ of length $A$ has no overlaps with any
relation in $\rho\setminus\rho_{T_i}$.

\section{Support varieties for simple modules}\label{sup_vari_for_sim}

In this section we give necessary and sufficient conditions for the variety
of a simple module to be nontrivial. 
We first recall the general fact that if 
$(P, \delta)$ is a minimal projective $\L^{e}$-resolution of $\L$, then,
for any simple $\L$-module $S$, we have that $(P\otimes_\L S,
\delta\otimes_\L S)$ is a minimal projective $\L$-resolution of $S$, and so
the complex $\Hom_\L(P\otimes_\L S,S)$ has zero differential.

\bigskip

For each vertex $v_j$ of the quiver ${\mathcal Q}$, for $1\leq j\leq
u$, let $S_j$ be the simple module corresponding to $v_j$. 
Let $A_j=\Ann_{\HH^*(\L)}\Ext^*_\L(S_j,S_j)$, so $A_j$ is a graded ideal
of $\HH^*(\L)$ and is properly contained in $\HH^*(\L)$.
Let $C_j$ be a closed path in ${\mathcal Q}$ at the vertex $v_j$, for some
$1 \leq j \leq u$. Then we say that $S_j$ is {\it associated to the 
closed path $C_j$}. We start by considering simple
modules associated to one of the closed paths $C_1, \ldots, C_u$.

\bigskip

\begin{lem}\label{ideal}
With the notation already introduced, fix an integer $j$ with $1\leq j\leq u$, so
that $C_j$ is a closed path in the quiver $\mathcal{Q}$ at the
vertex $v_j$. Let $z$ be a non-zero homogeneous element of
$\HH^*(\L)$ of the form $z=\sum_{i=1}^rc_ix_i^{l_i}$ with
$l_i\geq 0$, $c_i \in K$ and $c_j \neq 0$. Then $z$ is not
contained in $A_j$.
\end{lem}

\begin{proof}
If $\deg z=0$, then $z=c\cdot 1_{\HH^*(\L)}\in \HH^0(\L)\backslash
\operatorname{rad}\HH^0(\L)$ for some non-zero element $c\in K$, so
$z \not\in A_j$.

Now suppose $\deg z>0$, and assume for contradiction that $z\in
A_j$. Then $l_j>0$. Since $C_j$ is a closed path at the vertex
$v_j$, we know that $x_j$ is in degree $2$. Hence $\deg z=2l_j$.
Thus we have $\deg z = 2l_j =2l_t =2l_w \mu_w$ for all $1
\leq t \leq u$ and $u+1 \leq w \leq r$. We write $|z|$ for $\deg z$.
From the proof of \cite[Proposition 2.5]{GS}, for $1\leq t\leq u$,
the element $x_t^{l_t}\in \HH^{|z|}(\L)$ is represented by the
map $f_t: P^{|z|}\rightarrow\L$ where, for $R\in \R^{|z|}$,
\begin{equation}\label{eq1}
f_t\left(\mo(R)\otimes\mt(R)\right)=\left\{
\begin{array}{ll}
v_t & \mbox{if $R=C_t^{dl_t}$,}\\
0 & \mbox{otherwise,}
\end{array}
\right.
\end{equation}
and from the proof of \cite[Proposition 2.9]{GS}, for $u+1\leq w
\leq r$, the element $x_w^{l_w}\in \HH^{|z|}(\L)$ is represented
by the map $f_w: P^{|z|}\rightarrow\L$ where, for $R\in \R^{|z|}$,
\begin{equation}\label{eq2}
f_w\left(\mo(R)\otimes\mt(R)\right)=\left\{
\begin{array}{ll}
\mo(T_{w,k}) & \mbox{if $R=T_{w,k}^{dl_w/\gcd(d,m_w)}$,}\\
&\ \ \ \ \mbox{for $k = 0, \ldots , m_w-1$}\\
0 & \mbox{otherwise.}
\end{array}
\right.
\end{equation}
Then we have
$$f_i(\mo(C_j^{dl_j})\otimes \mt(C_j^{dl_j}))=\left\{
\begin{array}{ll}
v_j & \mbox{if $i=j$}\\
0 & \mbox{otherwise.}\\
\end{array}\right.
$$

Since $z \in A_j$, the element $z\otimes_\L S_j \in
\Ext_\L^{|z|}(S_j,S_j)$ is zero. Moreover since $z\otimes_\L
S_j$ is represented by $\sum_{i=1}^r (c_if_i\otimes_\L S_j):
P^{|z|}\otimes_\L S_j \rightarrow \L\otimes_\L S_j$, there is a left
$\L$-homomorphism $\Phi: P^{|z|-1}\otimes_\L S_j\rightarrow
\L\otimes_\L S_j$ such that $\sum_{i=1}^r (c_if_i\otimes_\L
S_j)=\Phi\circ(\delta^{|z|}\otimes_\L S_j)$.
So we have the following situation
$$\xymatrix{
P^{|z|}\otimes_\L S_j\ar[rr]^{\delta^{|z|}\otimes_\L
S_j}\ar[d]_{z\otimes_\L
S_j} & &  P^{|z|-1}\otimes_\L S_j\ar@{-->}[dll]^{\Phi}\\
\L\otimes_\L S_j & & \\
}$$
Then, by the remark at the start of this section, we have
$\Phi\circ(\delta^{|z|}\otimes_\L S_j)=0$. But
\begin{align*}
&\bigg(\sum_{i=1}^r (c_if_i\otimes_\L S_j)\bigg)(\mo(C_j^{dl_j})\otimes
\mt(C_j^{dl_j})\otimes_\L v_j)\\
&=(c_jf_j\otimes_\L S_j)(\mo(C_j^{dl_j})\otimes
\mt(C_j^{dl_j})\otimes_\L v_j) \\&= c_j(v_j\otimes_\L v_j)
\end{align*} which
is not zero in $\L\otimes_\L S_j$, since $c_j\neq0$. So
$\sum_{i=1}^r (c_if_i\otimes_\L S_j)$ is not zero. This is a
contradiction. Hence it follows that $z\not\in A_j$.
\end{proof}

\begin{prop}\label{closed_path}
With the notation already introduced, fix an integer $j$ with $1 \leq j\leq u$,
so that $C_j$ is a closed path in the quiver $\mathcal{Q}$ at the
vertex $v_j$. Then the variety of the simple module $S_j$ associated
to $C_j$ is nontrivial.
\end{prop}

\begin{proof}
We show first that $(A_j +\N)/\N$ is an ideal of the subalgebra
$$R_j=K[x_1, \ldots,x_{j-1},x_{j+1},\ldots , x_r]/\langle x_ax_b \mid a\neq b\rangle$$
of $K[x_1,\ldots,x_r]/\langle x_ax_b \mid a\neq b\rangle \cong
\HH^{*}(\L)/\N$.

Let $z+\N \in (A_j +\N)/\N$, so that $z = y + \chi$ for some $y \in A_j,
\chi \in \N$. Suppose first that $z$ is a homogeneous element of
$\HH^*(\L)$. Since $\HH^*(\L)/\N \cong K[x_1, \ldots , x_r]/\langle x_ax_b
\mid a \neq b\rangle$, we can also write
$z=\sum_{i=1}^r c_ix_i^{l_i} + \chi_1$ with $l_i\geq 0$, $c_i
\in K$ and $\chi_1 \in \N$. Thus $y = \sum_{i=1}^r c_ix_i^{l_i} + \chi_2$
where $\chi_2 = \chi_1 - \chi \in \N$. Since $\chi_2$ is nilpotent, there
is some $n \geq 0$ with $\chi_2^n = 0$, and so $(y - \sum_{i=1}^r
c_ix_i^{l_i})^n = 0$. Using the graded commutativity of $\HH^*(\L)$, we
then have that $(\sum_{i=1}^rc_ix_i^{l_i})^n \in \HH^*(\L)y$. Since $y$ is
in the ideal $A_j$ and $x_ax_b = 0$ for $a \neq b$, this gives
$\sum_{i=1}^rc_i^nx_i^{l_in} = (\sum_{i=1}^rc_ix_i^{l_i})^n \in A_j$.

Now, we may use Lemma \ref{ideal}, to show that either
$\sum_{i=1}^rc_i^nx_i^{l_in} = 0$ or $c_j^n = 0$. If
$\sum_{i=1}^rc_i^nx_i^{l_in} = 0$ then
$(\sum_{i=1}^rc_ix_i^{l_i})^n = 0$ so
$\sum_{i=1}^rc_ix_i^{l_i} \in \N$. Thus $z \in \N$ and so $z+\N =
0+\N \in R_j$. On the other hand, if $c_j^n = 0$ then $c_j = 0$ and
so $z+\N = \sum_{i=1,\ i\neq j}^rc_ix_i^{l_i} + \N \in R_j$. In
both cases we have that $z+\N \in R_j$.

Now suppose that $z$ is not a homogeneous element of $\HH^*(\L)$.
Since $A_j$ and $\N$ are both graded ideals, we have $z=\sum_{i=1}^Nz_i$
for some homogeneous elements $z_i = y_i + \chi_i$ with $y_i \in A_j,
\chi_i \in \N$ and some $N\geq 0$. From above, we know that
$z_i+\N \in R_j$ for $1\leq i\leq N$, and so $z+\N =
\sum_{i=1}^N(z_i+\N) \in R_j$.

Thus $(A_j +\N)/\N$ is contained in $R_j$. Therefore $(A_j+\N)/\N$
is contained in the maximal ideal $(x_1,\ldots x_{j-1},x_j-\alpha
,x_{j+1},\ldots,x_r)/ \langle x_ax_b\mid a\neq b\rangle$ of
$\HH^*(\L)/\N$ for all $\alpha \in K$.

Now,
\begin{align*}
V(S_j) &= \{\mm \in \MaxSpec\HH^*(\L)/\N \mid A_j \subseteq \mm'\}\\
&= \{\mm \in \MaxSpec\HH^*(\L)/\N \mid (A_j + \N)/\N \subseteq \mm\}.
\end{align*}
Thus $(x_1,\ldots x_{j-1},x_j-\alpha ,x_{j+1},\ldots,x_r)/
\langle x_ax_b\mid a\neq b\rangle \in V(S_j)$ for all
$\alpha \in K$. Hence the variety of $S_j$ is nontrivial.
\end{proof}

Let $T_j$ be a closed $A$-trail in the quiver ${\mathcal Q}$ for
some $u+1\leq j \leq r$. For $k$ with $0\leq k\leq m_{j-1}$, we denote by
$S_{jk}$ the simple module corresponding to the vertex
$\mo(T_{j,k})$. We say that the simple modules $S_{jk}$, for $0\leq
k\leq m_{j-1}$, are {\it associated to the closed $A$-trail $T_j$}.
Let $A_{jk}=\Ann_{\HH^*(\L)}\Ext^*_\L(S_{jk},S_{jk})$, so $A_{jk}$
is a graded ideal of $\HH^*(\L)$ and is properly contained in
$\HH^*(\L)$. We now show the analogue of Lemma~\ref{ideal} for
closed $A$-trails.

\begin{lem}\label{ideal_trail}
Keeping our notation, fix integers $j$ and $k$ with $u+1 \leq
j\leq r$ and $0\leq k\leq m_j-1$, so that $T_j$ is a closed
$A$-trail in the quiver $\mathcal{Q}$, and $T_{j,0}, \ldots ,T_{j,m_j-1}$
is the complete set of $A$-trails on $T_j$. Let $z$ be a non-zero
homogeneous element of $\HH^*(\L)$ of the form
$z=\sum_{i=1}^rc_ix_i^{l_i}$ with $l_i\geq 0$, $c_i \in K$ and
$c_j \neq 0$. Then $z$ is not contained in $A_{jk}$.
\end{lem}

\begin{proof} If $\deg z=0$, then $z=c\cdot 1_{\HH^*(\L)}\in \HH^0(\L)\backslash
\operatorname{rad}\HH^0(\L)$ for some non-zero element $c\in K$, so
$z \not\in A_{jk}$.

Now suppose that $\deg z>0$, and assume for contradiction that $z\in
A_{jk}$. Then $l_j>0$. Since $T_j$ is a closed $A$-trail, we know
that $x_j$ is in degree $2\mu_j$. Hence $\deg z=2l_j\mu_j$. Thus
we have $\deg z = 2l_j\mu_j=2l_t =2l_w \mu_w$ for all $1
\leq t \leq u$ and $u+1 \leq w \leq r$. Again we denote $\deg z$ by
$|z|$. Recall from the proof of Lemma~\ref{ideal} that, for $1\leq
t\leq u$, the element $x_t^{l_t}\in \HH^{|z|}(\L)$ is represented
by the map $f_t: P^{|z|}\rightarrow\L$ given by (\ref{eq1}), and,
for $u+1\leq w \leq r$, the element $x_w^{l_w}\in \HH^{|z|}(\L)$
is represented by the map $f_w: P^{|z|}\rightarrow\L$ given by
(\ref{eq2}). Then we have
$$f_i(\mo(T_{j,k}^{dl_j/\gcd(d,m_j)})\otimes \mo(T_{j,k}^{dl_j/\gcd(d,m_j)}))=\left\{
\begin{array}{ll}
\mo(T_{j,k}) & \mbox{if $i=j$}\\
0 & \mbox{otherwise.}\\
\end{array}\right.
$$

Since $z \in A_{jk}$, the element $z\otimes_\L S_{jk} \in
\Ext_\L^{|z|}(S_{jk},S_{jk})$ is zero. Moreover since $z\otimes_\L
S_{jk}$ is represented by $\sum_{i=1}^r (c_if_i\otimes_\L S_{jk}):
P^{|z|}\otimes_\L S_{jk} \rightarrow \L\otimes_\L S_{jk}$, there is
a left $\L$-homomorphism $\Phi: P^{|z|-1}\otimes_\L
S_{jk}\rightarrow \L\otimes_\L S_{jk}$ such that $\sum_{i=1}^r
(c_if_i\otimes_\L S_{jk})=\Phi\circ(\delta^{|z|}\otimes_\L S_{jk})$.
Then, by the remark at the start of this section, we have
$\Phi\circ(\delta^{|z|}\otimes_\L S_{jk})=0$. But
\begin{align*}
&\bigg(\sum_{i=1}^r (c_if_i\otimes_\L S_{jk})\bigg)(\mo(T_{j,k}^{dl_j/\gcd(d,m_j)})\otimes
\mo(T_{j,k}^{dl_j/\gcd(d,m_j)})\otimes_\L \mo(T_{j,k}))\\
&=(c_jf_j\otimes_\L S_{jk})(\mo(T_{j,k}^{dl_j/\gcd(d,m_j)})\otimes \mo(T_{j,k}^{dl_j/\gcd(d,m_j)})\otimes_\L \mo(T_{j,k}))\\
&=c_j(\mo(T_{j,k})\otimes_\L \mo(T_{j,k})),
\end{align*}
which is not zero in $\L\otimes_\L S_{jk}$, since $c_j\neq0$. So
$\sum_{i=1}^r (c_if_i\otimes_\L S_{jk})$ is not zero. This is a
contradiction. Hence it follows that $z\not\in A_{jk}$.
\end{proof}

Now, in a similar way to Proposition~\ref{closed_path}, we have the
following result concerning the varieties of simple modules associated to
closed $A$-trails.

\begin{prop}\label{closed_A_trail}
Keeping our notation, fix integers $j$ and $k$ with $u+1 \leq
j\leq r$ and $0\leq k\leq m_j-1$, so that $T_j$ is a closed
$A$-trail in the quiver $\mathcal{Q}$, and $T_{j,0}, \ldots ,T_{j,m_j-1}$
is the complete set of $A$-trails on $T_j$. Then the variety of the
simple module $S_{jk}$ is nontrivial.
\end{prop}

\begin{proof}
First we show that $(A_{jk} +\N)/\N$ is an ideal of the subalgebra
$$R_j=K[x_1, \ldots,x_{j-1},x_{j+1},\ldots , x_r]/\langle x_ax_b \mid
a \neq b\rangle$$ of $K[x_1,\ldots,x_r]/\langle x_ax_b\mid a \neq b
\rangle \cong \HH^{*}(\L)/\N$.

Let $z+\N \in (A_{jk} +\N)/\N$.
A similar argument to that used in the proof of Proposition~\ref{closed_path}
together with Lemma~\ref{ideal_trail} shows that $z + \N \in R_j$
Thus $(A_{jk} +\N)/\N$ is contained in $R_j$. Hence $(A_{jk}+\N)/\N$
is contained in the maximal ideal
$(x_1,\ldots x_{j-1},x_j-\alpha ,x_{j+1},\ldots,x_r)/
\langle x_ax_b\mid a \neq b\rangle$ of $\HH^*(\L)/\N$ for all $\alpha \in K$.

Again following Proposition ~\ref{closed_path}, we have that the
variety of $S_{jk}$ is nontrivial.
\end{proof}

Now we consider the varieties of modules which are not associated to
one of the closed paths $C_1,\ldots,C_u$ or the closed $A$-trails
$T_{u+1},\ldots,T_{r}$ in the quiver ${\mathcal Q}$.

\begin{prop}\label{not_on_closed_path}
Let $S$ be a simple module which is not associated to one of the
closed paths $C_1,\ldots,C_u$ or the closed $A$-trails
$T_{u+1},\ldots,T_{r}$ in the quiver ${\mathcal Q}$. Then the
variety of $S$ is trivial.
\end{prop}

\begin{proof}
If $\HH^*(\L)/\N \cong K$, then the variety of every finitely
generated module is trivial. In particular, the variety of each
simple module is trivial.

Assume now that $\HH^*(\L)/\N \not\cong K$, so we have some $r \geq 1$ with
$$\HH^*(\L)/\N \cong K[x_1, \ldots , x_r]/\langle x_ax_b
\mid a \neq b\rangle.$$ Recall that,
for $1\leq i \leq u$, the element $x_i$ corresponding to the closed
path $C_i$ at the vertex $v_i$ is represented by the map $P^2
\rightarrow \L$ where, for $R^2\in\R^2$,
$$\mo(R^2)\otimes \mt(R^2) \mapsto \left \{
\begin{array}{ll} v_i & \mbox{if $R^2 = C_i^d$}\\
                   0 & \mbox{otherwise,}
\end{array} \right.$$
and, for $u+1 \leq j \leq r$, the element $x_j$ corresponding to the
closed $A$-trail $T_j$ is represented by the map $P^{2\mu_j}
\rightarrow \L$ where, for $R^{2\mu_j}\in\R^{2\mu_j}$,
$$\mo(R^{2\mu_j})\otimes \mt(R^{2\mu_j}) \mapsto \left \{
\begin{array}{ll} \mo(T_{j,k}) & \mbox{if $R^{2\mu_j} = T_{j,k}^{d/\gcd(d,
m_j)}$}\\
&\ \ \ \ \mbox{for $k = 0, \ldots , m_j-1$}\\
                   0 & \mbox{otherwise.}
\end{array} \right.$$
Since $S$ is not associated to one of the closed paths $C_i$ at the
vertex $v_i$ for $1\leq i\leq u$ or to one of the closed $A$-trails
$T_j$ for $u+1\leq j\leq r$, it follows that $v_i \otimes_\L S=0$
for $1\leq i\leq u$, and $\mo(T_{j,k})\otimes_\L S=0$ for $u+1\leq
j\leq r$ and $0\leq k\leq m_j-1$. Thus the map $x_j\otimes_\L S$ is
zero for all $1\leq j\leq r$. So $(x_1,\ldots, x_r) \subseteq
\Ann_{\HH^*(\L)}\Ext_\L^*(S,S)$. Hence if $\mm \in
\MaxSpec\HH^*(\L)$ with $\Ann_{\HH^*(\L)}\Ext_\L^*(S,S)\subseteq
\mm$ then $\mm = (x_1,\ldots, x_r)$. Therefore $V(S)=\{(x_1, \ldots , x_r)\}$, 
so that the variety of $S$ is trivial.
\end{proof}

From Propositions~\ref{closed_path}, \ref{closed_A_trail}, and
\ref{not_on_closed_path}, we have the following theorem.

\begin{thm}\label{simple_module_variety}
Let $S$ be a simple module. Then the variety of $S$ is trivial if
and only if $S$ is not associated to one of the closed paths
$C_1,\ldots, C_u$ or to one of the closed $A$-trails
$T_{u+1},\ldots,T_r$ in the quiver ${\mathcal Q}$.
\end{thm}

We now give some examples to illustrate this theorem.

\begin{example}\label{example_1}
Let $\L = K{\mathcal Q}/I$ where ${\mathcal Q}$ is the quiver
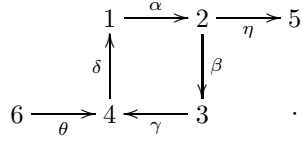
\[\xymatrix{%
& 1 \ar[r]^\alpha & 2 \ar[d]^\beta \ar[r]_\eta & 5 \\
6 \ar[r]_\theta& 4 \ar[u]^\delta & 3 \ar[l]^\gamma& \cdot}\] and $I$
is the ideal $\langle \alpha\beta, \beta\gamma, \gamma\delta,
\delta\alpha \rangle$. Then the simple modules $S_1$, $S_2$, $S_3$,
$S_4$ of $\L$ corresponding to the vertices $1$, $2$, $3$, $4$
respectively are associated to the closed $1$-trail
$\alpha\beta\gamma\delta$. So, by Theorem
\ref{simple_module_variety}, the varieties of $S_1$, $S_2$, $S_3$,
$S_4$ are nontrivial, whereas the varieties of the simple modules
corresponding to the vertices $5$, $6$ are trivial.
\end{example}

\begin{example}\label{example_2}
Let $\L = K{\mathcal Q}/I$ where ${\mathcal Q}$ is the quiver
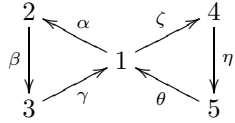
\[\xymatrix@R=6pt{%
2 \ar[dd]_\beta & & 4 \ar[dd]^\eta \\
  & 1 \ar[ul]_\alpha\ar[ur]^\zeta & \\
  3 \ar[ur]_\gamma & & 5 \ar[ul]^\theta}\]
and $I = \langle \alpha\beta, \beta\gamma, \gamma\alpha, \zeta\eta,
\eta\theta, \theta\zeta \rangle$. Then all simple modules of $\L$
are associated to closed $A$-trails in ${\mathcal Q}$. In fact, the
simple module corresponding to the vertex $1$ is associated to both
the closed $1$-trails $\alpha\beta\gamma$ and $\zeta\eta\theta$.
Thus, by Theorem~\ref{simple_module_variety}, the varieties of all
simple modules of $\L$ are nontrivial.
\end{example}

\begin{example}\label{example_3}
Let $\L = K{\mathcal Q}/I$ where ${\mathcal Q}$ is the quiver
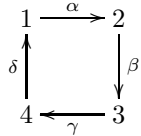
\[\xymatrix{ 1 \ar[r]^\alpha & 2 \ar[d]^\beta \\
4 \ar[u]^\delta & 3 \ar[l]^\gamma }\]
and $I
=\langle\alpha\beta\gamma\delta\alpha\beta,
\gamma\delta\alpha\beta\gamma\delta \rangle$. Then the simple
modules $S_1$, $S_3$ corresponding to the vertices $1$, $3$ are
associated to the closed $2$-trail $\alpha\beta\gamma\delta$. Hence,
by Theorem~\ref{simple_module_variety}, the varieties of $S_1$,
$S_3$ are nontrivial, whereas the varieties of the simple modules
corresponding to the vertices $2$, $4$ are trivial.
\end{example}

We complete this section with the following theorem which uses information
on the varieties of the simple modules to give structural information on
$\L$.

\begin{thm}\label{all_simples_nontrivial}
Let $\L=K\mathcal{Q}/I$ be a finite-dimensional $(D,A)$-stacked
monomial algebra. Suppose that each simple module has nontrivial variety.
Then $A=1$. Thus $\L$ is a $D$-Koszul monomial algebra.
\end{thm}

\begin{proof}
There is a simple $\L$-module with nontrivial variety, so
$\HH^*(\L)/\N \not\cong K$. Thus, from \cite[Theorem~3.4]{GS}, there
exists an integer $r\geq1$ with
$$\HH^*(\L)/\N \cong K[x_1, \ldots ,x_r]/\langle x_ax_b \mid a \neq b \rangle,$$
where $x_1$, \ldots , $x_u$ correspond to the closed trails $C_1$, \ldots , $C_u$ and
$x_{u+1}$, \ldots , $x_r$ correspond to the distinct closed $A$-trails
$T_{u+1}$, \ldots , $T_r$. From Theorem~\ref{simple_module_variety}, every simple
right $\L$-module is associated to a closed path $C_1$, \ldots , $C_u$ or a closed
$A$-trail $T_{u+1}$, \ldots , $T_r$. Thus, for each vertex $v$, we may fix a
relation $R^2_v$ in $\{C_1^d, \ldots , C_u^d\}\cup\left( \cup_{i=u+1}^{r} \rho_{T_i}\right)$
so that $R^2_v=vR^2_v$.

Choose a vertex $v_1\in \mathcal{Q}$. Then $v_1=\mo(R_{v_1}^2)$ for some
$R_{v_1}^2$ as above. Denote the first arrow in $R_{v_1}^2$ by $a_1$, and let
$v_2=\mt(a_1)$. Then $v_2=\mo(R_{v_2}^2)$. Again, let $a_2$ denote the first arrow
in $R_{v_2}^2$ and let $v_3=\mt(a_2)$. Continuing in this way gives a path
$a_1a_2a_3\cdots$ in $K\mathcal{Q}$. Since $\mathcal{Q}$ has a finite number
of vertices, eventually some vertex $v_i$ must repeat, so that there is a
closed path $a_1a_2a_3\cdots a_p$ in $K\mathcal{Q}$ with $\mo(a_1)=\mt(a_p)$
and $p\geq1$.

The algebra $\L$ is finite-dimensional, so there is a subpath of
$(a_1a_2\cdots a_p)^N$ for $N\gg0$ which lies in $\rho$, and so is
necessarily of length $D$. Denote this element of $\rho$ by $R_*^2$
and, by relabelling if necessary, suppose that the last arrow of
$R_*^2$ is $a_1$, so that $R_*^2=a_j\cdots a_{j+D-2}a_1$ for some
$j$. Thus we have the situation:
$$\xymatrix@W=0pt@M=0.3pt{
\ar@{^{|}-^{|}}@<-1.25ex>[rrr]_{R^2_*}\ar@{-}@<-1.25ex>[rr]^{a_j
\cdots a_{j+D-2}} & &
\ar@{-}@<-1.25ex>[r]^{a_1}\ar@{_{|}-_{|}}@<1.25ex>[rrr]^{R^2_{v_1}}
& & &}$$
Now $R_{v_1}^2$ overlaps $R_*^2$. But
$R_{v_1}^2\in\{C_1^d, \ldots ,C_u^d\}\cup\left(\cup_{i=u+1}^{r}
\rho_{T_i}\right)$.

If $R_{v_1}^2=C_i^d$ for some $1\leq i \leq u$, then $R_*^2=C_i^d$,
since we know that there are no overlaps of $C_i^d$ with any other
relation in $\rho$. But then, the first arrow of $R_*^2$ equals the
first arrow of $R_{v_1}^2$, so $a_j=a_1$ and $R_*^2=a_1a_{j+1}\cdots
a_{j+D-2}a_1=C_i^d$. But \cite[Lemma~2.2]{GS} gives us that
$\mo(R_*^2)=\mo(a_1)$ is not internal to $C_i$, so we must have
$R_*^2=a_1^D$ and $C_i=a_1$. Thus $A=\ell(C_i)=\ell(a_1)=1$.

Otherwise, we have $R_{v_1}^2\in \rho_{T_i}$ for some $u+1 \leq i
\leq r$. Since there are no overlaps of relations in $\rho_{T_i}$
with relations in $\rho\backslash \rho_{T_i}$, we know $R_*^2\in
\rho_{T_i}$ as well. By the construction of $\rho_{T_i}$, we know
that there is a relation $\tilde{R}^2_*$ in $\rho_{T_i}$ which ends
at $\mo(R_{v_1}^2)=v_1$. Then $R_*^2$ and $\tilde{R}^2_*$ lie in
$\rho_{T_i}$, and $R_*^2$ overlaps $\tilde{R}^2_*$ with overlap of
length $D+1$. This overlap is necessarily maximal and may be
illustrated as follows.
$$\xymatrix@W=0pt@M=0.3pt{
\ar@{^{|}-^{|}}@<-1.25ex>[rrr]_{\tilde{R}^2_*} &
\ar@{_{|}-_{|}}@<1.25ex>[rrr]^{R^2_*} & &
\ar@{^{|}-^{|}}@<-1.25ex>[rrr]_{R^2_{v_1}}\ar@{-}@<-1.25ex>[r]^{a_1}
& & & }$$
Hence there is an element of $\R^3$ of length $D+1$. Thus
$D+1=D+A$, so $A=1$.
\end{proof}

\section{Conditions {\bf (Fg1)} and {\bf (Fg2)}}\label{Section3}

In this section we give some examples of $(D,A)$-stacked monomial algebras 
which are not self-injective but nevertheless satisfy the following two
finiteness conditions found in \cite{EHSST}:
\begin{enumerate}
\item[\bf (Fg1)] There is a graded subalgebra $H$ of $\HH^*(\L)$ such that
\begin{enumerate}[(i)]
\item $H$ is a commutative Noetherian ring.
\item $H^0 = \HH^0(\L) = Z(\L)$.
\end{enumerate}
\item[\bf (Fg2)] The Ext algebra of $\L$
$$E(\L)=\Ext^*_\L(\L/\rad,\L/\rad)=\bigoplus_{i\geq 0}\Ext^i_\L(\L/\rad,\L/\rad)$$
is a finitely generated $H$-module.
\end{enumerate}
In the case where these conditions hold, we may apply the results of
\cite{EHSST} to give further significant information on support varieties
for $\L$-modules.
In particular, we can apply \cite[Theorem~2.5]{EHSST} to show that,
for a finite-dimensional $K$-algebra $\L$ satisfying {\bf (Fg1)}
and {\bf (Fg2)}, the algebra $\L$ is necessarily Gorenstein, and, for a finitely
generated $\L$-module $M$, the variety of $M$ is trivial if and
only if the projective dimension of $M$ is finite.

We now present some examples of $(D,A)$-stacked monomial algebras of
infinite global dimension which are not self-injective but which illustrate
the possible behaviours with respect to the conditions {\bf (Fg1)} and {\bf
(Fg2)}. In the first two examples, conditions {\bf (Fg1)} and {\bf
(Fg2)} hold.

\begin{example}
As in Example~\ref{example_2}, let $\L = K{\mathcal Q}/I$ where
${\mathcal Q}$ is the quiver
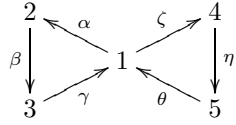
\[\xymatrix@R=6pt{%
2 \ar[dd]_\beta & & 4 \ar[dd]^\eta \\
  & 1 \ar[ul]_\alpha\ar[ur]^\zeta & \\
  3 \ar[ur]_\gamma & & 5 \ar[ul]^\theta}\]
and $I = \langle \alpha\beta, \beta\gamma, \gamma\alpha, \zeta\eta,
\eta\theta, \theta\zeta \rangle$. Then, $\L$ is a Koszul monomial
algebra (so is a $(2,1)$-stacked monomial algebra), and the variety
of every simple module is nontrivial.

Now, we claim that $\L$ satisfies {\bf (Fg1)} and {\bf (Fg2)}. By
\cite[Theorem 3.4]{GS}, $\HH^*(\L)/\N$ is isomorphic to the
subalgebra  $K[x,y]/(xy)$ of $\HH^*(\L)$, where $x$ corresponds to
the $1$-trail $\alpha\beta\gamma$ and is represented by the map $P^6
\rightarrow \L$ where, for
$R^6\in\R^6=\{(\alpha\beta\gamma)^2,(\beta\gamma\alpha)^2,
(\gamma\alpha\beta)^2,(\zeta\eta\theta)^2,(\eta\theta\zeta)^2,
(\theta\zeta\eta)^2\}$,
$$\mo(R^6) \otimes \mt(R^6)\mapsto
\left\{\begin{array}{ll}
e_1 & \mbox{if $R^6=(\alpha \beta \gamma)^2$},\\
e_2 & \mbox{if $R^6=(\beta \gamma \alpha)^2$},\\
e_3 & \mbox{if $R^6=(\gamma \alpha \beta)^2$},\\
0   & \mbox{otherwise},
\end{array}\right.$$
and $y$ corresponds to the $1$-trail $\zeta\eta\theta$ and is
represented by the map $P^6 \rightarrow \L$ where, for $R^6\in\R^6$,
$$
\mo(R^6) \otimes \mt(R^6) \mapsto \left\{\begin{array}{ll}
e_1 & \mbox{if $R^6=(\zeta \eta \theta)^2$},\\
e_4 & \mbox{if $R^6=(\eta \theta \zeta)^2$},\\
e_5 & \mbox{if $R^6=(\theta \zeta \eta)^2$},\\
0   & \mbox{otherwise}.
\end{array}\right.$$
It is easy to verify that $\HH^0(\L)=Z(\L)=K$. Let $H$ be the graded
subalgebra of $\HH^*(\L)$ generated by $1$, $x$ and $y$. Then it
follows that $H \cong K[x,y]/(xy)$ and $H^0= K$. Furthermore, $H$ is
a commutative Noetherian ring. Thus $\L$ satisfies {\bf (Fg1)}.

Now, we have $\R^0=\{e_1, \ldots ,e_5\}$ and $\R^1=\{\alpha, \beta,
\gamma, \zeta, \eta, \theta\}$. Since $\L$ is a Koszul monomial
algebra, we know that $E(\L)$ is generated by $\R^0$ and $\R^1$ as
an algebra. Let $n$ be a positive integer, and write $n=6q+r$ with
$0\leq r \leq 5$. For $i= 1$, $2$, $3$, let $P_{i,r}$ be the path of
length $r$ which lies on the closed $1$-trail $\alpha \beta \gamma$
such that $\mo(P_i)=i$. Similarly, for $i=1$, $4$, $5$, let
$Q_{i,r}$ be the path of length $r$ which lies on the closed
$1$-trail $\zeta \eta \theta$ such that $\mo(Q_i)=i$. Thus both
$P_{i,r}$ and $Q_{j,r}$ represent elements of degree $r$ in $E(\L)$
for $i=1$, $2$, $3$, and $j=1$, $4$, $5$. Moreover, since
$(\alpha\beta\gamma)^2$, $(\beta\gamma\alpha)^2$ and
$(\gamma\alpha\beta)^2$ represent elements of degree $6$ in $E(\L)$,
it follows that $(\alpha\beta\gamma)^{2q}$,
$(\beta\gamma\alpha)^{2q}$ and $(\gamma\alpha\beta)^{2q}$ represent
elements of degree $6q$ in $E(\L)$ for $n=6q+r$ with $0\leq r \leq
5$, where $(\alpha\beta\gamma)^0=e_1$, $(\beta\gamma\alpha)^0=e_2$
and $(\gamma\alpha\beta)^0=e_3$. By definition of $\R^n$, we get
$$\R^n=\{
(\alpha \beta \gamma)^{2q} P_{1,r}, (\beta \gamma \alpha)^{2q}
P_{2,r}, (\gamma \alpha \beta)^{2q} P_{3,r}, (\zeta \eta
\theta)^{2q} Q_{1,r}, (\eta \theta \zeta)^{2q}Q_{4,r}, (\theta \zeta
\eta)^{2q}Q_{5,r}\}$$ for $n=6q+r$ with $0\leq r\leq 5$. Also, we
know from Section 1 of \cite{GZ} that the set $\cup_{i\geq0}\R^i$ is
a multiplicative basis of $E(\L)$. Now the action of $H$ on $E(\L)$
is induced by the ring homomorphism $H \rightarrow E(\L)$ determined
by
$$\left\{\begin{array}{lll}
1 & \longmapsto & 1,\\
x & \longmapsto & (\alpha \beta \gamma)^2 + (\beta \gamma \alpha)^2
+(\gamma \alpha \beta)^2,\\
y & \longmapsto & (\zeta \eta \theta)^2 + (\eta \theta \zeta)^2 +
(\theta \zeta \eta)^2,
\end{array}\right.$$
so we have $\R^n = \{x^q P_{1,r}, x^q P_{2,r}, x^q P_{3,r}, y^q Q_{1,r},
y^q Q_{4,r}, y^qQ_{5,r} \}$ for $n=6q+r$ with $0 \leq r\leq 5$.
This shows that $E(\L)$ is generated by
$\cup_{0\leq i \leq 5} \R^i$ as a left $H$-module, and hence is a finitely
generated left $H$-module. Hence $\L$ satisfies {\bf (Fg2)}.

Moreover $e_1\L$ is not injective, so $\L$ is not self-injective. We
know from \cite[Theorem~2.5]{EHSST} that $\L$ is necessarily
Gorenstein, and it is easy to verify that the injective dimension of
$\L$ is $1$. Thus we have given an example of a nonself-injective
Koszul algebra where every simple module has nontrivial variety and 
{\bf (Fg1)} and {\bf (Fg2)} hold.
\end{example}

\begin{example}
As in Example~\ref{example_3}, let $\L = K{\mathcal Q}/I$ where
${\mathcal Q}$ is the quiver
\[\xymatrix{ 1 \ar[r]^\alpha & 2 \ar[d]^\beta \\
4 \ar[u]^\delta & 3 \ar[l]^\gamma }\]
and $I=\langle
\alpha\beta\gamma\delta\alpha\beta,
\gamma\delta\alpha\beta\gamma\delta \rangle$. Then the varieties of
the simple modules corresponding to $1$, $3$ are nontrivial, whereas
the varieties of the simple modules corresponding to $2$, $4$ are
trivial. Thus, for a simple module $S$, we see directly that the variety of
$S$ is trivial if and only if the projective dimension of $S$ is
finite. Moreover $e_3\L$ is not injective, so that $\L$ is not
self-injective.

We now show that $\L$ satisfies {\bf (Fg1)} and {\bf (Fg2)}, and
hence $\L$ must be Gorenstein by \cite[Theorem~2.5]{EHSST}. By
\cite[Theorem~3.4]{GS}, $\HH^*(\L)/\N \cong K[x]$ where $x$
corresponds to the 2-trail $\alpha \beta \gamma \delta$ and is
represented by the map $P^4 \rightarrow \L$, where, for $R^4 \in
\R^4 =\{(\alpha \beta \gamma \delta)^3, (\gamma \delta \alpha
\beta)^3\}$, $$\mo(R^4)\otimes \mt(R^4) \mapsto
\left\{\begin{array}{ll}
e_1 & \mbox{if $R^4=(\alpha \beta \gamma \delta)^3$}, \\
e_3 & \mbox{if $R^4=(\gamma \delta \alpha \beta)^3$}.
\end{array}\right.$$
Let $z=\alpha\beta\gamma\delta + \beta\gamma\delta\alpha +
\gamma\delta\alpha\beta + \delta\alpha\beta\gamma \in \L$. Then it
is easy to see that $z\in Z(\L)$ and $\HH^0(\L)=Z(\L)=K[z]/(z^2)$.
Moreover, let $H$ be the subalgebra of $\HH^*(\L)$ generated by $1$,
$z$, $x$. Then it follows that $H^0 = K[z]/(z^2)$ and $H \cong
K[x,z]/(z^2)$, and hence $H$ is a commutative Noetherian ring.
Therefore $\L$ satisfies {\bf (Fg1)}.

Now, we have $\R^0=\{e_1, e_2, e_3, e_4\}$ and $\R^1=\{\alpha, \beta,
\gamma, \delta\}$. Let
$$R_1=\alpha\beta\gamma\delta\alpha\beta, \quad
R_2=\gamma\delta\alpha\beta\gamma\delta, \quad
R'_1=\alpha\beta\gamma\delta\alpha\beta\gamma\delta, \quad
R'_2=\gamma\delta\alpha\beta\gamma\delta\alpha\beta.$$
Then, for $i=1$, $2$, $R_i$ and $R'_i$ represent elements of degrees
$2$ and $3$ of $E(\L)$ respectively, and we have $\R^2=\{R_1, R_2 \}$ and
$\R^3=\{R'_1, R'_2 \}$.  We know from \cite{GS2} that $E(\L)$ is generated by
$\R^0$, $\R^1$, $\R^2$, and $\R^3$ as an algebra.
Let $n$ be an integer with $n\geq4$, and write $n=4q+r$ with
$0\leq r \leq 3$ (so $q\geq1$). Then, we get
$$\R^n = \left\{\begin{array}{ll}
\{(R_1R_2)^q, (R_2R_1)^q\} & \mbox{if $r=0$},\\
\{(R_1R_2)^{q-1}R_1R'_2, (R_2R_1)^{q-1}R_2R'_1\}& \mbox{if $r=1$},\\
\{(R_1R_2)^qR_1,(R_2R_1)^qR_2\}& \mbox{if $r=2$},\\
\{(R_1R_2)^qR'_1,(R_2R_1)^qR'_2\}& \mbox{if $r=3$}.
\end{array}\right.
$$
Also, we know from Section 1 of \cite{GZ} that the set $\cup_{i\geq0}\R^i$
is a multiplicative basis of $E(\L)$. Now the action of $H$ on $E(\L)$
is induced by the ring homomorphism $H\rightarrow E(\L)$ determined by
$$\left\{\begin{array}{lll}
1 & \longmapsto & 1,\\
x & \longmapsto & R_1R_2 + R_2R_1,\\
z & \longmapsto & 0,
\end{array}\right.$$
so that, for $n\geq4$, we have
$$\R^n = \left\{\begin{array}{ll}
\{x^q e_1, x^q e_2\} & \mbox{if $r=0$},\\
\{x^{q-1}R_1R'_2, x^{q-1}R_2R'_1\}& \mbox{if $r=1$},\\
\{x^qR_1, x^qR_2\}& \mbox{if $r=2$},\\
\{x^qR'_1, x^qR'_2\}& \mbox{if $r=3$}.
\end{array}\right.$$
This shows that $E(\L)$ is  generated by
$\cup_{0\leq i \leq 5} \R^i$ as a left $H$-module. 
So $\L$ satisfies {\bf (Fg2)}.

Moreover it can be seen that the injective dimension of $\L$ is $2$.
In fact, $e_2\L$ and $e_4\L$ are injective modules, and moreover we
have the following injective resolutions of $e_1\L$ and $e_3\L$
respectively:
$$0\rightarrow e_1 \L
\longrightarrow I(2) \longrightarrow I(4) \longrightarrow I(3)
\rightarrow 0,$$
$$0\rightarrow e_3 \L \longrightarrow I(4) \longrightarrow I(2) \longrightarrow
I(1) \rightarrow 0,$$ where $I(j)$ is the indecomposable injective
module corresponding to the vertex $j$ for $1 \leq j \leq 4$. Thus
the injective dimensions of $e_1\L$ and $e_3\L$ are $2$, and so the
injective dimension of $\L$ is $2$. Thus we have given an example of
a nonself-injective monomial algebra where some simple modules have
trivial variety, some simple modules have nontrivial variety, and
{\bf (Fg1)} and {\bf (Fg2)} hold.
\end{example}

In the next two examples, conditions {\bf (Fg1)} and {\bf
(Fg2)} do not hold.

\begin{example}
Let ${\mathcal Q}$ be the quiver
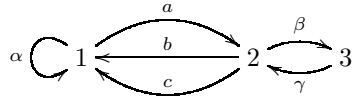
\[\xymatrix{
1 \ar@/^5mm/[rr]^a \ar@(ul,dl)[]_\alpha && 2 \ar[ll]_b
\ar@/^5mm/[ll]_c \ar@/^/[r]^\beta & 3 \ar@/^/[l]^\gamma }\] and let
$\L = K{\mathcal Q}/I$ where $I =\langle \alpha^2, \beta\gamma,
\gamma\beta, ab, ac, ba, ca \rangle$. Then $\L$ is Koszul (so is
a $(2,1)$-stacked monomial algebra). Also, the vertex $1$ is
associated to the closed path $\alpha$, and the vertices $2$ and $3$
are associated to the $1$-trail $\beta \gamma$. So, by
Theorem~\ref{simple_module_variety}, the variety of each simple
module is nontrivial.

Now, there is the following minimal projective resolution of the indecomposable
injective right $\L$-module $I(3)$ corresponding to the vertex $3$:
$$\cdots \longrightarrow U_3 \longrightarrow U_2 \longrightarrow U_1
\longrightarrow e_3\L \oplus e_3\L \longrightarrow I(3) \rightarrow0,$$
where, for $i\geq1$,
$$U_i=
\left\{\begin{array}{ll}
(e_1\L)^{\oplus 3\cdot 2^{m-1}} & \mbox{if $i=2m-1$}, \\
(e_2\L)^{\oplus 3\cdot 2^{m-1}} & \mbox{if $i=2m$}.
\end{array}\right.$$
Hence, the projective dimension of $I(3)$ is infinite, which implies
that the injective dimension of $\L e_3$ is infinite. So $\L$ is not
Gorenstein. Hence, by \cite[Theorem~2.5]{EHSST}, there is no $H$ satisfying
{\bf (Fg1)} and {\bf (Fg2)}. Therefore we have given an example of a
nonself-injective Koszul monomial algebra where the variety of every
simple module is nontrivial, and {\bf (Fg1)} and {\bf (Fg2)} do not
hold.
\end{example}

\begin{example}
Let ${\mathcal Q}$ be the quiver
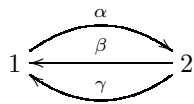
\[\xymatrix{
1 \ar@/^5mm/[rr]^\alpha  && 2 \ar[ll]_\beta \ar@/^5mm/[ll]_\gamma
}\]
and let $\L = K{\mathcal Q}/I$ where $I =\langle \alpha \beta,
\alpha \gamma, \gamma \alpha, \beta \alpha \rangle$. Then $\L$ is a
Koszul monomial algebra with radical squared zero. From
\cite[Theorem~3.4]{GS}, $\HH^*(\L)/\N \cong K$, and so the variety of every
finitely generated module is trivial. Thus all simple modules of
$\L$ have trivial variety. Moreover, it is easy to see that all the
simple $\L$-modules have infinite projective dimension. It now
follows from \cite[Theorem~2.5]{EHSST} that there is no $H$ such
that {\bf (Fg1)} and {\bf (Fg2)} hold. In particular, $\L$ is not
Gorenstein. Thus we have given an example of a nonself-injective
Koszul monomial algebra where the variety of every simple module is
trivial, and {\bf (Fg1)} and {\bf (Fg2)} do not hold.
\end{example}

\end{document}